\documentclass[12pt,english]{article}
\usepackage[T1]{fontenc}
\usepackage[latin1]{inputenc}
\usepackage[a4paper]{geometry}
\usepackage{amsmath}
\usepackage{amssymb}
\usepackage{graphicx}
\usepackage{esint}

\makeatletter
\@ifundefined{date}{}{\date{}}
\usepackage{amsthm}\usepackage[compress]{cite}
\usepackage{cleveref}

\theoremstyle{plain}

\theoremstyle{plain}

\theoremstyle{plain}

\DeclareMathOperator{\re}{Re}
\DeclareMathOperator{\im}{Im}

\providecommand{\lemmaname}{Lemma}
\providecommand{\propositionname}{Proposition}
\providecommand{\corollaryname}{Corollary}

\makeatother

\usepackage{babel}
\usepackage{url}

\begin{document}
\title{\textbf{Rayleigh waves in isotropic strongly elliptic thermoelastic
materials with microtemperatures}}
\author{F. Passarella\thanks{Dipartimento di Matematica, Università di Salerno, Italy - email:
fpassarella@unisa.it.} \and {}V. Tibullo\thanks{Dipartimento di Matematica, Università di Salerno, Italy - email:
vtibullo@unisa.it.}\thanks{Corresponding author.} \and {}G. Viccione\thanks{Dipartimento di Ingegneria Civile, Università di Salerno, Italy -
email: gviccion@unisa.it.}}
\maketitle
\begin{abstract}
This paper is concerned with the linear theory of thermoelasticity
with microtemperatures, based on the entropy balance proposed by Green
and Naghdi, which permits the transmission of heat as thermal waves
of finite speed. We analyze the behavior of Rayleigh waves in an unbounded
isotropic homogeneous strongly elliptic thermoelastic material with
microtemperatures. The related solution of the Rayleigh surface wave
problem is expressed as a linear combination of the elements of the
bases of the kernels of appropriate matrices. The secular equation
is established and afterwards an explicit form is written when some
coupling constitutive coefficients vanish. Then, we solve numerically
the secular equation by mean of a grafical metod and by taking arbitrary
data for strongly elliptic thermoelastic material\emph{.}

Published in \\
Passarella, F., Tibullo, V. \& Viccione, G. \emph{Rayleigh waves in isotropic strongly elliptic thermoelastic materials with microtemperatures}. Meccanica \textbf{52}, 3033--3041 (2017). \\
\url{https://doi.org/10.1007/s11012-016-0591-z}
\end{abstract}
\textbf{\small{}Keywords}{\small{}: Rayleigh waves; strong ellipticity;
microtemperatures.}{\small\par}

\section{Introduction}

There are many authors who study the materials having thermal variations
at the microstructural level, such as viscous fluids, granular materials,
composites and nanomaterials (e.g. \cite{Eringen1999,Iesan2004,jVadaszGovenderpVadasz2005,JordanPuri2001}
ans references). This is due to increasing interest in several class
of nanomaterials used in the heat transfer industry, where the microtemperatures
and microdeformations of the nanoparticles cannot be ignored. In future
technologies studies related to propagation wave in the theory of
thermoelastic materials with microtemperatures may be important.

On the other hand, during the past years several authors have studied
the class of strongly elliptic materials. These materials are characterized
by special properties, like negative Poisson's ratio and negative
stiffness (auxetic or antirubber materials). These particular structures
(see for example \cite{ParkLakes2007}) expand laterally when stretched,
in contrast to the behavior of ordinary materials. The ellipticity
analysis is relevant in studying wave propagation \cite{Gurtin1972}
and has important applications in several contexts (e.g. \cite{MerodioOgden2003,MerodioOgden2005,CiarlettaChiritaPassarella2005,TibulloVaccaro2008,PassarellaZampoli2009a,PassarellaZampoli2009b,PassarellaTibulloZampoli2010,PassarellaTibulloZampoli2011}).

The propagation of thermoelastic waves has been discussed long ago
by Lockett \cite{Lockett1958} and Lockett and Sneddon \cite{LockettSneddon159}.
Chadwick \cite{Chadwick1960} studied the coupled and modified character
of the thermoelastic waves and noted that they are also damped. Later,
Ivanov \cite{Ivanov1988} used these results to discuss appropriate
criteria for the behavior at infinity, in order to preserve the characteristic
features of the Rayleigh waves known from the classical elasticity. 

The effects of heat conduction upon the propagation of Rayleigh surface
waves in a semi-infinite elastic solid has been studied by Chadwick
and Windle\cite{ChadwickWindle1964} in isotropic thermoelastic bodies
and by Chakraborty and Pal \cite{ChakrabortyPal1969} and by Chadwick
and Seet \cite{Chadwick Seet1970} for transversely isotropic materials.
Further, Abouelregal\cite{Abouelregal} studied Rayleigh waves in
a thermoelastic homogeneous solid half space in the context of a dual-phase-lag
model. We have to point out that the wave motion in the form of acceleration
waves and of shock waves is discussed in the recent book by Straughan
\cite{Straughan2011} in an account of theories of heat conduction
where the temperature may travel as a wave with finite speed. Further,
the propagation of elastic waves and the propagation of Rayleigh surface
waves has been studied in various contexts, in \cite{Deresiewicz1957,Achenbach1967,Puri1972,Agarwal1979,IovaneNasedkinPassarella2004,IovaneNasedkinPassarella2005,ChiritaCiarlettaTibullo2013,BucurPassarellaTibullo2014,CiarlettaPassarellaSvanadze2014,ChiritaCiarlettaTibullo}. 

In \cite{CiarlettaPassarellaTibullo2015}, Ciarletta \emph{et al.}
study a homogeneous strongly elliptic thermoelastic body with microtemperatures
following the theory of Iesan and Quintanilla \cite{IesanQuintanilla2009}.
In \cite{IesanQuintanilla2009} it is presented a linearized theory
based on the entropy balance proposed by Green and Naghdi \cite{GreenNaghdi1991}.
Moreover, in \cite{CiarlettaPassarellaTibullo2015} the authors show
that there is neither dispersion nor attenuation in the wave propagation
as a consequence of the entropy balance proposed by Green and Naghdi
\cite{GreenNaghdi1991,GreenNaghdi1993} and the strong ellipticity
condition; this is in contrast to what we see in \cite{SteebSinghTomar2012}
where the theory of thermoelasticity with microtemperatures of Iesan
and Quintanilla \cite{IesanQuintanilla2000} is used. Further, the
authors prove that only undamped plane harmonic waves exist for any
direction of propagation. In the isotropic case the possible waves
are undamped in time and there are three longitudinal and two transverse
waves. We point out that all transverse waves have constant temperature.

In the present paper, the theory of thermoelasticity with microtemperatures
(Iesan and Quintanilla \cite{IesanQuintanilla2009}, Ciarletta et
al. \cite{CiarlettaPassarellaTibullo2015}) is applied to the study
of Rayleigh waves propagating at the thermally insulated stress-free
surface of an isotropic, homogeneous strongly elliptic thermoelastic
solid half-space with microtemperatures.

The layout of the paper is a follows. In Section 2, we state the set
of basic equations describing the behavior of thermoelastic media
with microtemperatures within the context of the theory developed
in \cite{IesanQuintanilla2009}. Further, we remark some of the results
obtained in\cite{CiarlettaPassarellaTibullo2015} as the conditions
characterizing the strong ellipticity for isotropic materials.

In Section 3, we find the explicit solutions for surface waves propagation
in a half space filled with an isotropic homogeneous strongly elliptic
thermoelastic medium with microtemperatures. We prove that there is
no dispersion, moreover we obtain that the solution of the Rayleigh
surface wave problem is expressed as a linear combination of the elements
of the (five) bases of the kernels of the appropriate matrices; these
vectors are written in an explicit form.\textbf{ }The secular equation
is then established, then we solve numerically by mean of a grafical
metod and by taking arbitrary data for strongly elliptic thermoelastic
material.

Finally, in Section 4 all vectors of the bases of the considered kernels,
and the corresponding secular equations, are established in three
cases in which some of the coupling constitutive coefficients vanish. 

\section{Field equations}

Let $\Omega$ be an unbounded region filled of a thermoelastic material
with microstructure, as presented in \cite{IesanQuintanilla2009,CiarlettaPassarellaTibullo2015}.
For this type of bodies, the temperature $\theta^{\prime}$ at the
point $\mathbf{X}^{\prime}$ of the microelement $\omega$ is considered
a linear function of the microcoordinates $\mathbf{X}^{\prime}-\mathbf{X}$,
i.e. $\theta^{\prime}=\theta+\mathbf{T\cdot}(\mathbf{X}^{\prime}-\mathbf{X})$,
where $\mathbf{X}$ is the center of mass of $\omega$ in the reference
configuration, $\theta$ and $\mathbf{T}$ are the temperature and
the microtemperature vector at $\mathbf{X}$ and $\cdot$ is the euclidean
scalar product. It is assumed that there exists a reference time $t_{R}$
such that
\[
\mathbf{T}(\mathbf{X},t_{R})=\mathbf{\mathbf{T}}^{R},\qquad\theta(\mathbf{X},t_{R})=\theta^{R}.
\]

In the following, a rectangular Cartesian coordinate system $Ox_{k}$,
$k=1,2,3$, is used. Letters in boldface, like $\mathbf{v}$, stand
for tensors of order $p$, with components $v_{ij...s}$ ($p$ subscripts).
Latin subscripts range over the integers $\{1,2,3\}$, Greek subscripts
range over $\{1,2\}$ and the summation convention is employed. A
superposed dot or a subscript preceded by a comma will mean partial
derivative with respect to time or to the corresponding coordinate,
respectively. Moreover, we suppress the dependence upon spatial and/or
temporal variables when no confusion may occur. All involved functions
are supposed to be sufficiently regular to ensure analysis to be valid.

In the context of the linear theory presented in \cite{IesanQuintanilla2009,CiarlettaPassarellaTibullo2015}
and in the absence of supply terms, the behavior of the isotropic
homogeneous body possessing a center of symmetry is governed by the
following equations 
\begin{equation}
\begin{alignedat}{1} & \mu u_{i,jj}+\left(\lambda+\mu\right)u_{j,ji}+\varepsilon_{2}\tau_{i,jj}+\left(\varepsilon_{1}+\varepsilon_{2}\right)\tau_{j,ji}-\beta\dot{\chi}_{,i}=\rho\ddot{u}_{i},\\
 & \varepsilon_{2}u_{i,jj}+\left(\varepsilon_{1}+\varepsilon_{2}\right)u_{j,ji}+d_{2}\tau_{i,jj}+\left(d_{1}+d_{3}\right)\tau_{j,ji}-m\dot{\chi}_{,i}=b\ddot{\tau}_{i},\\
 & -\beta\dot{u}_{j,j}-m\dot{\tau}_{j,j}+k\chi_{,jj}=a\ddot{\chi}.
\end{alignedat}
\qquad\text{on }\Omega\times(t_{R},\infty),\label{eq:isotropic-motion-eq}
\end{equation}
where ${\bf u}=(u_{1},\,u_{2},\,u_{3})^{T}$ is the displacement vector
field and $\boldsymbol{\tau}=(\tau_{1},\,\tau_{2},\,\tau_{3})^{T}$
and $\chi$ are defined by 
\[
\tau_{i}=\int_{t_{R}}^{t}\left[T_{i}(s)-T_{i}^{R}\right]ds,\qquad\chi=\int_{t_{R}}^{t}\left[\theta(s)-\theta^{R}\right]ds.
\]
Here $\rho$ is the reference mass density and $\lambda$, $\mu$,
$\varepsilon_{1}$, $\varepsilon_{2}$, $d_{1}$, $d_{2}$, $d_{3}$,
$\beta$, $m$, $k$, $a$ and $b$ are constitutive coefficients. 

We remark that, for the considered theory, the coupling constitutive
coefficients are $\varepsilon_{1}$, $\varepsilon_{2}$, $\beta$
and $m$. 

We can see that the equations of the system \eqref{eq:isotropic-motion-eq}
are uncoupled when

\global\long\def\labelenumi{\roman{enumi})}%

\begin{enumerate}
\item $\beta=0$, $m\neq0$ and $\varepsilon_{1}=\varepsilon_{2}=0$ (Eq.
\eqref{eq:isotropic-motion-eq}$_{1}$ reduces to the classical motion
equation of elasticity);
\item $\beta\neq0$, $m=0$ and $\varepsilon_{1}=\varepsilon_{2}=0$ (Eq.
\eqref{eq:isotropic-motion-eq}$_{2}$ involves only $\boldsymbol{\tau}$);
\item $\beta=0$, $m=0$, $\varepsilon_{1}\neq0$ and $\varepsilon_{2}\neq0$
(Eq. \eqref{eq:isotropic-motion-eq}$_{3}$ involves only $\chi$).
\end{enumerate}
It is proved in \cite{CiarlettaPassarellaTibullo2015} that a thermoelastic
material with microtemperatures is strongly elliptic if and only if
\begin{equation}
\rho>0,\qquad a>0,\qquad b>0,\qquad k>0,\label{eq:strongly-elliptic-1}
\end{equation}
\begin{equation}
\lambda+2\mu>0,\qquad\mu>0,\qquad(\varepsilon_{1}+2\varepsilon_{2})^{2}<(\lambda+2\mu)d,\qquad\varepsilon_{2}^{2}<\mu d_{2},\label{eq:strongly-elliptic-2}
\end{equation}
where $d=d_{1}+d_{2}+d_{3}.$ Consequently, it is also $d>0$, $d_{2}>0$. 

In what follows, it is useful to introduce the following polynomials
\begin{equation}
q_{2}(t)=t^{2}-a_{2}t+a_{0},\qquad q_{3}(t)=t^{3}-b_{4}t^{2}+b_{2}t-b_{0},\qquad t\in\mathbb{C}\label{eq:polynomials}
\end{equation}
with
\[
\begin{alignedat}{1} & a_{2}=\frac{\mu}{\rho}+\frac{d_{2}}{b},\qquad a_{0}=\frac{\mu d_{2}-\varepsilon_{2}^{2}}{\rho b},\\
 & b_{4}=\left(\dfrac{\lambda+2\mu}{\rho}+\dfrac{d}{b}\right)+\dfrac{1}{a}\left(\dfrac{m^{2}}{b}+\dfrac{\beta^{2}}{\rho}\right)+\dfrac{k}{a},\\
 & b_{2}=\dfrac{1}{\rho abd}\left\{ \left(ad+m^{2}\right)\left[\left(\lambda+2\mu\right)d-\left(\varepsilon_{1}+2\varepsilon_{2}\right)^{2}\right]+\left[d\beta-(\varepsilon_{1}+2\varepsilon_{2})m\right]^{2}\right\} \\
 & \quad+\dfrac{k}{a}\left(\dfrac{\lambda+2\mu}{\rho}+\dfrac{d}{b}\right),\\
 & b_{0}=\dfrac{k}{\rho ab}\left[(\lambda+2\mu)d-\left(\varepsilon_{1}+2\varepsilon_{2}\right)^{2}\right].
\end{alignedat}
\]
 Under the restrictions imposed by the strong ellipticity conditions
\eqref{eq:strongly-elliptic-1}, \eqref{eq:strongly-elliptic-2},
we can easily prove that 
\[
a_{2}>0,\qquad a_{0}>0,\qquad b_{4}>0,\qquad b_{2}>0,\qquad b_{0}>0,
\]
so that if the polynomials have a real root, it must be positive.
Moreover, it is proved in \cite{CiarlettaPassarellaTibullo2015} that
if the constitutive coefficients satisfy the conditions \eqref{eq:strongly-elliptic-1}
and \eqref{eq:strongly-elliptic-2}, then the equation 
\[
q_{2}(t)q_{3}(t)=0
\]
has only real (and positive) solutions. In particular, we have that
the roots of $q_{2}(t)$ are
\begin{equation}
\begin{alignedat}{1} & t_{1,2}=\dfrac{\mu b+\rho d_{2}\pm\sqrt{\left(\mu b-\rho d_{2}\right)^{2}+4\rho b\varepsilon_{2}^{2}}}{2\rho b}.\end{alignedat}
\label{eq:roots-q2}
\end{equation}
On the other hand, in \cite{CiarlettaPassarellaTibullo2015} it is
shown that the cubic equation $q_{3}(t)=0$ has three different solutions
if 
\begin{equation}
h_{0}^{2}<\dfrac{4}{27}h_{1}^{3},\label{eq:6}
\end{equation}
with $h_{0}=-(2b_{4}^{3}-9b_{2}b_{4}+27b_{0})/27,$ and $h_{1}=(b_{4}^{2}-3b_{2})/3,$
and these roots are 
\begin{equation}
\begin{alignedat}{1} & t_{k}=\dfrac{b_{4}}{3}+2\sqrt{\frac{h_{1}}{3}}\cos\left[\dfrac{1}{3}\arccos\left(\frac{-3h_{0}}{2h_{1}}\sqrt{\frac{3}{h_{1}}}\right)-\frac{2\pi}{3}(k+1)\right]\end{alignedat}
\quad k=3,4,5.\label{eq:roots-q3}
\end{equation}
In the following, we assume that the constitutive coefficients are
such that Eqs. \eqref{eq:strongly-elliptic-1}, \eqref{eq:strongly-elliptic-2}
and \eqref{eq:6} are satisfied.

\section{Rayleigh surface waves}

In what follows, $\Omega$ is a half-space made of an isotropic homogeneous
strongly elliptic thermoelastic material with microtemperatures such
that the coupling constitutive coefficients are non zero ($m\neq0$,
$\beta\neq0$ and $\varepsilon_{1}\neq0$ or $\varepsilon_{2}\neq0$).
This half-space is characterized by $x_{2}\geq0$. We study the propagation
of a surface wave in the $x_{1}$-direction and with attenuation in
the $x_{2}$-direction. The surface $x_{2}=0$ is assumed stress free
and thermally insulated. To this end, we seek for solutions of the
system \eqref{eq:isotropic-motion-eq} in the form 
\begin{equation}
\begin{alignedat}{3}u_{\alpha} & =u_{\alpha}\left(x_{1}-vt,x_{2}\right), & \quad u_{3} & =0,\\
\tau_{\alpha} & =\tau_{\alpha}\left(x_{1}-vt,x_{2}\right), & \tau_{3} & =0, & \quad\chi & =\chi\left(x_{1}-vt,x_{2}\right),\quad\text{with }v\in\mathbb{C},
\end{alignedat}
\label{eq:8}
\end{equation}
so that the system (\ref{eq:isotropic-motion-eq}) becomes 
\begin{equation}
\begin{alignedat}{1} & \left[t_{1i}-\rho v^{2}u_{i,1}\right]_{,1}+\left[t_{2i}\right]_{,2}=0,\\
 & \left[\Lambda_{1i}-v^{2}b\tau_{i,1}\right]_{,1}+\left[\Lambda_{2i}\right]_{,2}=0,\\
 & \left[S_{1}-v^{2}a\chi_{,1}\right]_{,1}+\left[S_{2}\right]_{,2}=0,
\end{alignedat}
\label{eq:9}
\end{equation}
where 
\begin{equation}
\begin{alignedat}{1}\begin{alignedat}{1} & t_{\alpha i}=\lambda\delta_{i\alpha}u_{\beta,\beta}+\mu(u_{\alpha,i}+u_{i,\alpha})+\varepsilon_{1}\delta_{i\alpha}\tau_{\beta,\beta}+\varepsilon_{2}(\tau_{\alpha,i}+\tau_{i,\alpha})+v\beta\delta_{i\alpha}\chi_{,1},\\
 & \Lambda_{\alpha i}=\varepsilon_{1}\delta_{\alpha i}u_{\beta,\beta}+\varepsilon_{2}(u_{\alpha,i}+u_{i,\alpha})+d_{1}\delta_{i\alpha}\tau_{\beta,\beta}+d_{2}\tau_{i,\alpha}+d_{3}\tau_{\alpha,i}+vm\delta_{i\alpha}\chi_{,1},\\
 & S_{\alpha}=v\beta u_{\alpha,1}+vm\tau_{\alpha,1}+k\chi_{,\alpha}.
\end{alignedat}
\end{alignedat}
\label{eq:constitutiveEq}
\end{equation}
Let be 
\begin{equation}
\begin{aligned}{\cal U} & =(u_{1},\,u_{2},\,\tau_{1},\,\tau_{2},\,\chi)^{T}, & \quad{\cal T}_{\alpha} & =(t_{\alpha1},\,t_{\alpha2},\,\Lambda_{\alpha1},\,\Lambda_{\alpha2},\,S_{\alpha})^{T}.\end{aligned}
\label{eq:11}
\end{equation}
Since we study the wave with attenuation in the direction $x_{2}$
and the surface $x_{2}=0$ is stress free and is thermally insulated,
we have the following asymptotic conditions 
\begin{equation}
\lim_{x_{2}\to+\infty}{\cal U}\left(x_{1},x_{2},t\right)=0,\qquad\lim_{x_{2}\to+\infty}{\cal T}_{\alpha}\left(x_{1},x_{2},t\right)=0,\qquad\forall x_{1}\in\mathbb{R},\,t\geq0,\label{eq:asymptotic-cond.}
\end{equation}
and the boundary conditions
\begin{equation}
{\cal T}_{2}\left(x_{1},0,t\right)=0,\qquad\forall x_{1}\in\mathbb{R},\,t\geq0.\label{eq:boubdarycond}
\end{equation}
Now, we seek for solutions ${\cal U}$ of the above problem in the
following exponential form
\begin{equation}
{\cal U}=\tilde{{\cal U}}e^{i\varkappa\left(x_{1}-vt+px_{2}\right)}\text{\qquad with }\tilde{{\cal U}}=(U_{1},\,U{}_{2},\,A_{1},\,A_{2},\,B)^{T}.\label{eq:exp-form}
\end{equation}
Here, ${\bf U}=(U_{1},\,U{}_{2},\,0)^{T}$ and $\mathbf{A}=(A_{1},\,A_{2},\,0)^{T}$
are complex constant vectors and $B$ is a complex constant with $\left\vert \mathbf{U}\right\vert \neq0$
or $\left\vert {\bf A}\right\vert \neq0$ or $B\neq0$. Further, $v$
is such that 
\begin{equation}
v=v_{R}-v_{I}i\qquad\text{with }v_{R}\geq0,\,v_{I}\geq0,\label{eq:ReIm}
\end{equation}
and $p$ is such that 
\begin{equation}
p=\alpha+\beta i\qquad\text{with }\beta>0,\label{eq:p}
\end{equation}
in order to satisfy the asymptotic conditions \eqref{eq:asymptotic-cond.}.
In particular, the real part of $v$ gives the wave speed and the
imaginary part gives the rate of damping in time.

Substituting the exponential form \eqref{eq:exp-form} into Eqs. \eqref{eq:9},
\eqref{eq:constitutiveEq} we arrive to a homogeneous linear algebraic
system
\begin{equation}
\mathfrak{D}_{p}\tilde{{\cal U}}=0,\label{eq:linear-system}
\end{equation}
where the matrix $\mathfrak{D}_{p}$ is defined as 
\begin{equation}
\mathfrak{D}_{p}=p^{2}{\cal Q}_{1}+p{\cal Q}_{2}+{\cal R}\label{eq:matrixD_p}
\end{equation}
with
\[
\!\begin{alignedat}{1}{\cal Q}_{1} & =\left(\begin{array}{ccccc}
\mu & 0 & \varepsilon_{2} & 0 & 0\\
0 & \lambda+2\mu & 0 & \varepsilon_{1}+2\varepsilon_{2} & 0\\
\varepsilon_{2} & 0 & d_{2} & 0 & 0\\
0 & \varepsilon_{1}+2\varepsilon_{2} & 0 & d & 0\\
0 & 0 & 0 & 0 & k
\end{array}\right),{\cal \;Q}_{2}=\left(\begin{array}{ccccc}
0 & \lambda+\mu & 0 & \varepsilon_{1}+\varepsilon_{2} & 0\\
\lambda+\mu & 0 & \varepsilon_{1}+\varepsilon_{2} & 0 & v\beta\\
0 & \varepsilon_{1}+\varepsilon_{2} & 0 & d_{1}+d_{3} & 0\\
\varepsilon_{1}+\varepsilon_{2} & 0 & d_{1}+d_{3} & 0 & mv\\
0 & v\beta & 0 & mv & 0
\end{array}\right),\\
\\
{\cal R} & =\left(\begin{array}{ccccc}
\lambda+2\mu-\rho v^{2} & 0 & \varepsilon_{1}+2\varepsilon_{2} & 0 & v\beta\\
0 & \mu-v^{2}\rho & 0 & \varepsilon_{2} & 0\\
\varepsilon_{1}+2\varepsilon_{2} & 0 & d_{0}-bv^{2} & 0 & mv\\
0 & \varepsilon_{2} & 0 & d_{2}-bv^{2} & 0\\
v\beta & 0 & mv & 0 & k-av^{2}
\end{array}\right),
\end{alignedat}
\]
in other words we should have $\tilde{{\cal U}}\in\ker\mathfrak{D}_{p}$.
Since ${\cal \tilde{{\cal U}}}$ is a non-trivial solution of the
homogeneous algebraic system \eqref{eq:linear-system}, then we obtain
the following propagation condition
\begin{equation}
\det\mathfrak{D}_{p}=0.\label{eq:detD_p}
\end{equation}
We observe that $\det\mathfrak{D}_{p}$ is a homogeneous polynomial
with respect to $v^{2}$ and $p^{2}+1$; consequently, $v=0$ if and
only if $p^{2}+1=0$. Now, if we introduce 
\begin{equation}
t=\frac{v^{2}}{p^{2}+1}\qquad\text{with}\;v\neq0,\label{eq:t}
\end{equation}
we can rewrite the condition \eqref{eq:detD_p} as 
\begin{equation}
q_{2}(t)q_{3}(t)=0,\label{eq:detD_p-1}
\end{equation}
where $q_{2}$ and $q_{3}$ are defined in \eqref{eq:polynomials}. 

The propagation condition \eqref{eq:detD_p-1} holds if and only if
$t$ is a root of $q_{2}$ or $q_{3}$. We suppose that $q_{2}$ and
$q_{3}$ do not have a common root. The (real and positive) solutions
of Eq. \eqref{eq:detD_p-1} are expressed in Eqs. \eqref{eq:roots-q2}
and \eqref{eq:roots-q3}. 

The wave-number $\varkappa$ does not appear in Eq. \eqref{eq:detD_p},
so that the phase velocity $v$ cannot depend on $\varkappa$ and
therefore there is no dispersion. 

Let be $t_{k}$ a solution of Eq. \eqref{eq:detD_p-1} ($k=1,\ldots,5$)
and be $p_{k}=\alpha_{k}+\beta_{k}i$ the value of $p$ corresponding
to $t_{k}$ through Eq. \eqref{eq:t} and satisfying Eq. \emph{\eqref{eq:p}}.
Then, Eqs. \emph{\eqref{eq:ReIm}} and \eqref{eq:t} imply
\begin{equation}
t_{k}\left[\alpha_{k}^{2}-\beta_{k}^{2}+1\right]=v_{R}^{2}-v_{I}^{2},\qquad t{}_{k}\alpha_{k}\beta_{k}=-v_{R}v_{I},\quad\forall k.\label{eq:22}
\end{equation}
Since all roots $t_{k}$ are positive, we arrive to 
\begin{equation}
\alpha_{k}\beta_{k}\le0\qquad\implies\alpha_{k}\le0.\label{eq:pk}
\end{equation}
Consequently, we can remark by using Eqs.\eqref{eq:22} that 
\[
v_{R}v_{I}=0\qquad\implies\alpha_{k}=0\quad\text{or\quad}\beta_{k}=0;
\]
in particular, we obtain from \eqref{eq:22}$_{1}$:
\begin{description}
\item [{$v=v_{R}:$}] $p_{k}=\sqrt{1-\dfrac{v_{R}^{2}}{t_{k}}}\,i$ if
$v_{R}<\min_{k\in\left\{ 1,\ldots,5\right\} }\sqrt{t_{k}}$, and $p_{k}=\alpha_{k}$
if $v_{R}\geq\sqrt{t_{k}}\text{ for some }k$ but it is not compatible
with the condition \eqref{eq:p};
\item [{$v=v_{I}i:$}] $p_{k}=\beta_{k}i\quad\forall k$.
\end{description}
The solution \eqref{eq:exp-form} of the problem corresponding to
$p_{k}$ $(k=1,\ldots,5)$ is
\begin{equation}
{\cal U}^{(k)}={\cal \tilde{U}}^{(k)}e^{i\varkappa\left(x_{1}-vt+p_{k}x_{2}\right)}={\cal {\cal \tilde{U}}}^{(k)}e^{-\varkappa\beta_{k}x_{2}}e^{i\varkappa\left(x_{1}-vt-\alpha_{k}x_{2}\right)},\label{eq:exp-form-1}
\end{equation}
where 
\begin{equation}
{\cal \tilde{U}}{}^{(k)}=(U{}_{1}^{(k)},\,U{}_{2}^{(k)},\,A_{1}^{(k)},\,A{}_{2}^{(k)},\,B^{(k)})^{T}\in\ker\mathfrak{D}_{p_{k}}.\label{eq:kerD_pk}
\end{equation}
We calculate the solutions of the corresponding homogeneous linear
system \eqref{eq:linear-system} and we arrive to 
\[
\begin{alignedat}{2}\tilde{{\cal U}}{}^{(k)}= & \left(-p_{k}\Phi_{k},\,\Phi_{k},\,-p_{k},\,1,\,0\right)^{T}, &  & \text{\quad if }k=1,2,\\
\\
\tilde{{\cal U}}{}^{(k)}= & \Biggl\{\Gamma_{k},\,p_{k}\Gamma_{k},\,\Lambda_{k},\,p_{k}\Lambda_{k},\,\dfrac{v}{m\beta t_{k}}\left[\Gamma_{k}\Lambda_{k}-\left(\varepsilon_{1}+2\varepsilon_{2}\right)(\beta\Gamma_{k}+m\Lambda_{k})\right]\Biggr\}^{T}, &  & \text{\quad\ if }k=3,4,5,
\end{alignedat}
\]
with 
\[
\begin{array}{ll}
\Phi_{k}=\dfrac{b}{\varepsilon_{2}}\left(t_{k}-\dfrac{d_{2}}{b}\right),\\
\\
\Gamma_{k}=b\beta\left(t_{k}-\dfrac{d_{2}}{b}\right)+m\left(\varepsilon_{1}+2\varepsilon_{2}\right), & \;\Lambda_{k}=\rho m\left(t_{k}-\dfrac{\lambda+2\mu}{\rho}\right)+\mathbf{\beta}\left(\varepsilon_{1}+2\varepsilon_{2}\right).
\end{array}
\]
Let be $\mathbf{n}{}^{(k)}$$=\left(1,\,p_{k},\,0\right)^{T}$. It
is then obvious that:
\begin{description}
\item [{$k=1,\:2$:}] $\tilde{{\cal U}}{}^{(k)}$ are such that $\mathbf{U}{}^{(k)}=(-p_{k}\Phi_{k},\,\Phi_{k},\,0)^{T}$
and ${\bf A}{}^{(k)}$$=(-p_{k},\,1,\,0)$ are orthogonal to $\mathbf{n}{}^{(k)}$;
\item [{$k=3,\:4,\:5$:}] $\tilde{{\cal U}}{}^{(k)}$are such that ${\bf U}{}^{(k)}$$=(\Gamma_{k},\,p_{k}\Gamma_{k},\,0)^{T}$
and ${\bf A}{}^{(k)}=(\Lambda_{k},\,p_{k}\Lambda_{k},\,0)$ are parallel
to $\mathbf{n}{}^{(k)}$.
\end{description}
The more general solution ${\cal U}$ of our problem is given by a
linear combination of the ${\cal U}^{(k)}$
\begin{equation}
{\cal U}(x_{1},x_{2},t)=\sum_{k=1}^{5}\gamma_{k}\tilde{{\cal U}}^{(k)}e^{i\varkappa(x_{1}-vt+p_{k}x_{2})},\label{eq:linear-comb}
\end{equation}
where $\boldsymbol{\gamma}=(\gamma_{1},\,\gamma_{2},\,\ldots,\,\gamma_{5})^{T}$
is a non-zero constant vector. Substituting the expression \eqref{eq:linear-comb}
into Eq. \eqref{eq:constitutiveEq}, we arrive to
\begin{equation}
{\cal T}_{2}(x_{1},x_{2},t)=i\varkappa\sum_{k=1}^{5}\gamma_{k}{\cal S}_{p_{k}}\tilde{{\cal U}}^{(k)}e^{i\varkappa\left(x_{1}-vt+p_{k}x_{2}\right)}\label{eq:T_2}
\end{equation}
with

\[
{\cal S}_{p_{k}}=\left(\begin{array}{ccccc}
\mu p_{k} & \mu & p_{k}\varepsilon_{2} & \varepsilon_{2} & 0\\
\lambda & (\lambda+2\mu)p_{k} & \varepsilon_{1} & p_{k}\left(\varepsilon_{1}+2\varepsilon_{2}\right) & v\beta\\
p_{k}\varepsilon_{2} & \varepsilon_{1} & d_{2}p_{k} & d_{3} & 0\\
\varepsilon_{2} & p_{k}\left(\varepsilon_{1}+2\varepsilon_{2}\right) & d_{1} & dp_{k} & mv\\
0 & v\beta & 0 & mv & kp_{k}
\end{array}\right).
\]
With the aid of the boundary conditions \eqref{eq:boubdarycond},
Eq. \eqref{eq:T_2} leads to
\[
\mathcal{T}_{2}(x_{1},0,t)=i\varkappa\sum_{k=1}^{5}\gamma_{k}{\cal S}_{p_{k}}\tilde{{\cal U}}{}^{(k)}e^{i\varkappa\left(x_{1}-vt\right)}=0,\qquad\forall x_{1}\in\mathbb{R},\forall t\geq0,
\]
 and, equivalently,
\begin{equation}
{\cal A}{\bf \gamma}=0\label{eq:Agamma}
\end{equation}
where ${\cal A}=\left\Vert a_{hk}\right\Vert $ with $(a_{1k},a_{2k},...,a_{5k})^{T}={\cal S}_{p_{k}}{\cal U}{}^{(k)}$.
A non trivial solution $\boldsymbol{\gamma}$ of Eq. \eqref{eq:Agamma}
exists if and only if 
\begin{equation}
\det{\cal A}=0,\label{eq:detAgamma}
\end{equation}
which represents the secular equation for the complex parameter $v$.
We have to select the solutions of the secular equation \eqref{eq:detAgamma}
satisfying the conditions \eqref{eq:ReIm}.

Now, we want to investigate, from a numerical point of view, the secular
equation \eqref{eq:detAgamma} with respect to the complex parameter
$v$. To this aim, we will take arbitrary values for the relevant
constitutive parameters, compatible with restrictions \eqref{eq:strongly-elliptic-1},
\eqref{eq:strongly-elliptic-2}. In particular, we look for a numerical
solution of Eq. \eqref{eq:detAgamma}.

We note that this relation contains the unknown $v$ both explicitly
and implicitly through $p_{k}$, $k=1,\ldots,5$, that should be taken
as the solutions of relation \eqref{eq:detD_p}. The solution of this
system of two nonlinear equations is not easy, and we take another
approach. We define 
\[
\mathcal{F}(\re(v),\im(v))=\ln|\det{\cal A}|.
\]
The presence of the logarithm is convenient because the function has
a wide range of variability. We can now make a graphics of the function
$\mathcal{F}$ looking for a minimum. In Figure 1 we show the graphics
that we have obtained, where it is possible to see the presence of
a minimum around $v=0.62-0.08i$.

\begin{center}
\begin{figure}
\centering\includegraphics[scale=0.7]{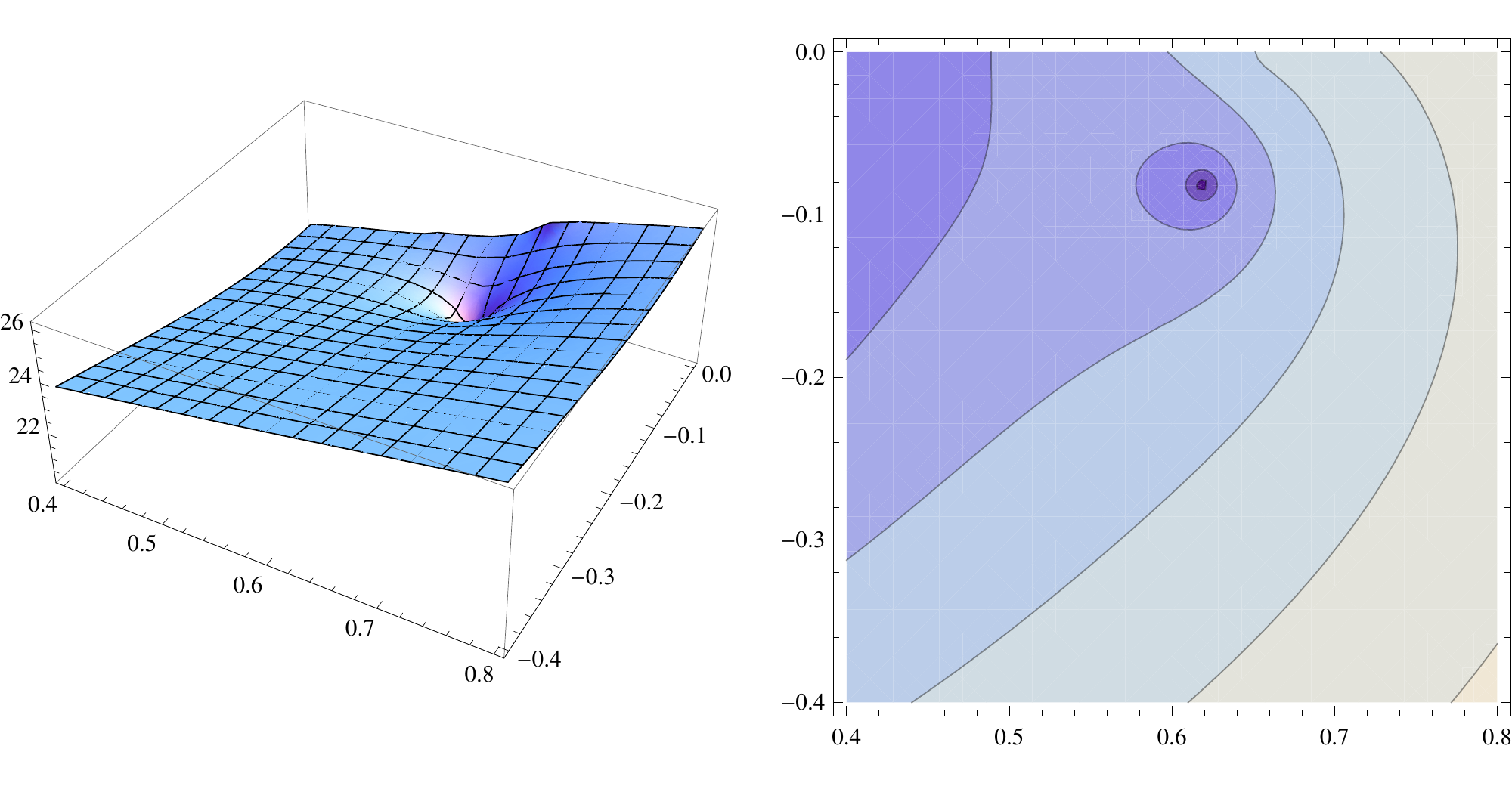}\caption{The graphics of the function $\mathcal{F}(\re(v),\im(v))$ for $\re(v)\in(0.4,\,0.8)$
and for $\im(v)\in(-0.4,\,0.0)$.}
\end{figure}
\par\end{center}

\section{Special cases }

In this section, we consider the class of isotropic strongly elliptic
thermoelastic media with microtemperatures when some of the coupling
coefficients vanish; in particular, we consider the following cases:

\global\long\def\labelenumi{\roman{enumi})}%

\begin{enumerate}
\item $\beta=0$, $m\neq0$ and $\varepsilon_{1}=\varepsilon_{2}=0$:

The propagation condition \eqref{eq:detD_p} is rewritten as Eq. \eqref{eq:detD_p-1},
where $q_{2}$ and $q_{3}$ reduce to
\begin{equation}
q_{2}=\left(t-\frac{\mu}{\rho}\right)\left(t-\frac{d_{2}}{b}\right),\;q_{3}=\left(t-\frac{\lambda+2\mu}{\rho}\right)\left(t^{2}-\frac{m^{2}+ad+bk}{ab}t+\frac{kd}{ab}\right),\label{eq:q2q3caso1}
\end{equation}
so that 
\begin{equation}
t_{1}=\frac{\mu}{\rho},\;t_{2}=\frac{d_{2}}{b},\label{eq:t-12-caso1}
\end{equation}
 
\begin{equation}
t_{3}=\frac{\lambda+2\mu}{\rho},\;t_{4,5}=\frac{m^{2}+ad+bk\pm\sqrt{m^{4}+\left(ad-bk\right)^{2}+2m^{2}\left(ad+bk\right)}}{2ab}.\label{eq:t-345-caso1}
\end{equation}
Let be $p_{k}$ the values of $p$ corresponding, through \eqref{eq:t},
to the roots $t_{k}$ defined in \eqref{eq:t-12-caso1} and \eqref{eq:t-345-caso1}.
The kernels, associated with $\mathfrak{D}_{p_{k}}$, are spanned
by 
\[
\begin{alignedat}{1}\tilde{{\cal U}}^{(1)} & =\left(-p_{1},\,1,\,0,\,0,\,0\right)^{T},\\
\tilde{{\cal U}}^{(2)} & =\left(0,\,0,\,-p_{2},\,1,\,0\right),^{T}\\
\tilde{{\cal U}}^{(3)} & =\left(1,\,p_{3},\,0,\,0,\,0\right)^{T},\\
\tilde{{\cal U}}^{(4)} & =\left(0,\,0,\,\varPi_{4},\,p_{4}\varPi_{4},\,mv(bt_{4}-d_{2})\right)^{T},\\
\tilde{{\cal U}}^{(5)} & =\left(0,\,0,\,\varPi_{5},\,p_{5}\varPi_{5},\,mv(bt_{5}-d_{2})\right)^{T},
\end{alignedat}
\]
where 
\[
\varPi_{k}=m^{2}t_{k}+\left(at_{k}-k\right)(d_{1}+d_{3}),\qquad k=4,5.
\]

The secular equation \eqref{eq:detAgamma} reduces to

\[
\begin{aligned} & v\left[4\mu^{2}p_{1}p_{2}+\left(\rho v^{2}-2\mu\right)^{2}\right]\Bigl\{ bp_{4}\left[bv^{2}-\left(d_{2}+d_{3}\right)\right]\left[\left(k-at_{5}\right)\left[bv^{2}-\left(d_{2}+d_{3}\right)\right]+m^{2}v^{2}\right]\\
 & +p_{5}\biggl[p_{3}p_{4}\left(d_{2}+d_{3}\right){}^{2}\left(-2abt_{5}+ad+bk\right)+a\left(d-bt_{5}\right)\left[bv^{2}-\left(d_{2}+d_{3}\right)\right]{}^{2}\\
 & +\left.m^{2}\left(d_{2}+d_{3}\right)\left[\left(p_{3}p_{4}+1\right)\left(d_{2}+d_{3}\right)-bv^{2}\right]\right]\Bigr\}=0.
\end{aligned}
\]

\item $\beta\neq0$, $m=0$ and $\varepsilon_{1}=\varepsilon_{2}=0$:

The propagation condition \eqref{eq:detD_p} leads to Eq. \eqref{eq:detD_p-1}
where $q_{2}$ is defined by Eq. \eqref{eq:q2q3caso1} and $q_{3}$
reduces to 
\[
q_{3}(t)=\left(bt-d\right)\left(t^{2}-\frac{a(\lambda+2\mu)+\beta^{2}+k\rho}{a\rho}t+\frac{(\lambda+2\mu)k}{a\rho}\right).
\]
In particular, the roots $t_{1},t_{2}$ are defined by \eqref{eq:t-12-caso1}
and the other roots are
\[
\begin{aligned}t_{3}= & \frac{d}{b},\qquad t_{4,5}=\frac{a(\lambda+2\mu)+\beta^{2}+k\rho\pm\sqrt{\beta^{4}+\left[a(\lambda+2\mu)-\rho k\right]^{2}+2\beta^{2}\left[a(\lambda+2\mu)+\rho k\right]}}{2a\rho}.\end{aligned}
\]
Consequently, we obtain 
\[
\begin{alignedat}{1}\tilde{{\cal U}}^{(1)} & =\left(-p_{1},\,1,\,0,\,0,\,0\right)^{T},\\
\tilde{{\cal U}}^{(2)} & =\left(0,\,0,\,-p_{2},\,1,\,0\right)^{T},\\
\tilde{{\cal U}}^{(3)} & =\left(0,\,0,\,1,\,p_{3},\,0\right),^{T}\\
\tilde{{\cal U}}^{(4)} & =\left(\varOmega_{4},\,p_{4}\varOmega_{4},\,0,\,0,\,\beta v(\rho t_{4}-\mu)\right)^{T},\\
\tilde{{\cal U}}^{(5)} & =\left(\varOmega_{5},\,p_{5}\varOmega_{5},\,0,\,0,\,\beta v(\rho t_{5}-\mu)\right),^{T}
\end{alignedat}
\]
where 
\[
\varOmega_{k}=\beta^{2}t_{k}+\left(at_{k}-k\right)(\lambda+\mu),\qquad k=4,5.
\]
We can calculate the secular equation \eqref{eq:detAgamma} and we
obtain
\[
\begin{aligned} & v\left[\left[bv^{2}-\left(d_{2}+d_{3}\right)\right]{}^{2}+p_{2}p_{3}\left(d_{2}+d_{3}\right){}^{2}\right]\Bigl\{ p_{4}\rho\left(\rho v^{2}-2\mu\right)\left[\beta v^{2}-\left(k-at_{5}\right)\left(2\mu-\rho v^{2}\right)\right]\\
 & +p_{5}\Bigl[4\mu^{2}p_{1}p_{4}\left(a\left(\lambda+2\mu-2\rho t_{5}\right)+k\rho\right)+a\left(\rho v^{2}-2\mu\right)^{2}\left(\lambda+2\mu-\rho t_{5}\right)\\
 & +2\beta^{2}\mu\left(2\mu+2\mu p_{1}p_{4}-\rho v^{2}\right)\Bigr]\Bigr\}=0.
\end{aligned}
\]

\item $\beta=0$, $m=0$, $\varepsilon_{1}\neq0$ and $\varepsilon_{2}\neq0$:

The propagation condition \eqref{eq:detD_p} leads to Eq. \eqref{eq:detD_p-1},
where $q_{2}$ is defined by Eq. \eqref{eq:q2q3caso1} and $q_{3}$
reduces to 
\[
q_{3}(t)=\left(t-\frac{k}{a}\right)\left(t^{2}-\dfrac{b\left(\lambda+2\mu\right)+d\rho}{\rho b}t+\dfrac{(\lambda+2\mu)d-\left(\varepsilon_{1}+2\epsilon_{2}\right)^{2}}{\rho b}\right).
\]
In particular, the roots $t_{1},t_{2}$ are defined by \eqref{eq:roots-q2}
and other three roots are 
\[
t_{3}=\dfrac{k}{a},\qquad t_{4,5}=\frac{b\left(\lambda+2\mu\right)+d\rho\pm\sqrt{\left[b(\lambda+2\mu)-\rho d\right]^{2}+4\rho b\left(\varepsilon_{1}+2\varepsilon_{2}\right)^{2}}}{2\rho b}.
\]
The vectors of the bases of the kernels of $\mathfrak{D}_{p_{k}}$
are defined by 
\[
\begin{alignedat}{1}\tilde{{\cal U}}^{(1)} & =\left(\varepsilon_{1}+2\varepsilon_{2},\,(\varepsilon_{1}+2\varepsilon_{2})p_{1},\,\varPsi_{1},\,p_{2}\varPsi_{1},\,0\right)^{T},\\
\tilde{{\cal U}}^{(2)} & =\left(\varepsilon_{1}+2\varepsilon_{2},\,(\varepsilon_{1}+2\varepsilon_{2})p_{2},\,\varPsi_{2},\,p_{2}\varPsi_{2},\,0\right)^{T},\\
\tilde{{\cal U}}^{(3)} & =\left(-\varepsilon_{2}p_{3},\,\varepsilon_{2},\,p_{3}\hat{\varPsi}_{3},\,\hat{\varPsi}_{3},\,0\right)^{T},\\
\tilde{{\cal U}}^{(4)} & =\left(-\varepsilon_{2}p_{4},\,\varepsilon_{2},\,p_{4}\hat{\varPsi}_{4},\,\hat{\varPsi}_{4},\,0\right)^{T},\\
\tilde{{\cal U}}^{(5)} & =\left(0,\,0,\,0,\,0,\,\varepsilon_{2}\right)^{T},
\end{alignedat}
\]
where 
\[
\begin{alignedat}{2}\varPsi_{k} & =\rho t_{k}-(\lambda+2\mu),\quad &  & k=1,2\\
\hat{\varPsi}_{k} & =\rho t_{k}-\mu &  & k=3,4
\end{alignedat}
\]
It is possible to calculate the secular equation \eqref{eq:detAgamma},
but the obtained formula is too long to be reported here.
\end{enumerate}

\end{document}